\documentclass[12pt,a4paper,reqno]{amsart}
\usepackage{amsmath,amsthm}
\usepackage{amsfonts}
\usepackage{amssymb}
\allowdisplaybreaks

\newcommand{\be}{\begin{equation}}
\newcommand{\ef}{\end{equation}}
\chardef\bslash=`\\ 

\newtheorem*{thm*}{Theorem}

\theoremstyle{definition}

\newtheorem*{remark*}{Remarks}
\newtheorem*{defn*}{Definition}

\theoremstyle{remark}

\newcommand{\wt}{\widetilde}
\newcommand{\wh}{\widehat}
\newcommand{\fc}{\frac}

\newcommand{\iy}{\infty}
 \renewcommand{\sectionmark}[1]{}

\newcommand{\ve}{\varepsilon}

\newcommand{\const}{\operatorname{const}}

 \usepackage{amsfonts}
\newcommand{\field}[1]{\mathbb{#1}}

\newcommand{\D}{\field{D}}
\newcommand{\om}{\omega}
\newcommand{\z}{\zeta}
\newcommand{\ov}{\overline}
\newcommand{\vp}{\varphi}
\newcommand{\hC}{\wh{\field{C}}}
\newcommand{\C}{\field{C}}

\newcommand{\B}{\mathbf{B}}
\newcommand{\T}{\mathbf{T}}
\newcommand{\Belt}{\operatorname{Belt}}

\newcommand{\grad}{\operatorname{grad}}

\newcommand{\x} {\mathbf x}
\renewcommand{\a} {\alpha}

\begin{document}
\title{Teichm\"{u}ller balls and biunivalent holomorphic functions}
\author{Samuel L. Krushkal}

\begin{abstract} Biunivalent holomorphic functions form an interesting class in geometric function theory
and are connected with special functions and solutions of complex differential equations.

The paper reveals a deep connection between biunivalence and geometry of Teichm\"{u}ller balls
and provides some sufficient conditions for biunivalence of holomorphic functions on the disk. Among the
consequences, one obtains new sharp distortion theorems for univalent functions with quasiconformal
extension.
\end{abstract}

\date{\today\hskip4mm({Tbiuniv.tex})}

\maketitle

\bigskip

{\small {\textbf {2020 Mathematics Subject Classification:} Primary: 30C25, 30C45, 30C50, 30C62,
30G70, 30F60; Secondary 30C62, 32F45}

\medskip

}

\bigskip

\markboth{Samuel L. Krushkal}{Teichm\"{u}ller balls and biunivalent holomorphic functions} \pagestyle{headings}

\bigskip\noindent
{\bf 1. Preamble}. Biunivalent functions form an interesting subclass of holomorphic functions
in the unit disk $\D = \{|z| < 1\}$ which consists of univalent functions $w = f(z) = z + a_2 z^2 + \dots$
in $\D$ whose inverse functions $z = f^{-1}(w)$ also are univalent in this disk. Such functions have a closed  connection with special functions, solutions of complex differential equations, etc. and still are investigated from these points of view. Biunivalence causes rather strong rigidity.

The aim of this paper is to describe geometrically the place of these functions among all univalent functions.
A characterization is given using the features of geodesic balls $B_\kappa$ in the universal
Teichm\"{u}ller space.

The result is that all points of the balls of radii $\kappa \le 1/4$ represent biunivalent
functions. The upper bound $1/4$ is sharp and cannot be increased.

This shows that the collection of such functions is rather broad and, on the other hand, complements
the known results on univalence and quasiconfornmal extensions of holomorphic functions, which play a fundamental role in complex analysis (see \cite{AW}, \cite{Be}, \cite{Ber}, \cite{N}).
Other applications of our main theorem concern the extremal problems for univalent functions with
wuasiconformal extension.

\bigskip\noindent
{\bf 2. A glimpse to universal Teichm\"{u}ller space}.
Recall that the  universal Teichm\"{u}ller space $\T =\T(\D)$ is the space of quasisymmetric homeomorphisms of the unit circle $\mathbb S^1 = \partial \D$ factorized by M\"{o}bius maps.
Its canonical complex Banach structure is defined by factorization of the ball of Beltrami coefficients
$$
\Belt(\D)_1 = \{\mu \in L_\iy(\C) : \ \mu|\D^* = 0, \ \| \mu \| < 1\},
$$
letting $\mu, \nu \in \Belt(\D)_1$ be equivalent if the corresponding quasiconformal maps $w^{\mu_1},
w^{\mu_2}$ (defined as the solutions of the Beltrami equation $\partial_{\ov z} w = \mu(z) \partial_z w$
with the hydrodynamic normalization $w(\iy) = \iy, \ w^\prime(\iy) = 1$ for $\mu = \mu_1, \ \mu_2$)
coincide on $\mathbb S^1$.
The {\bf dilatation} $k(w) = \|\mu\|_\iy$ estimates the deviation of $w^\mu(z)$ from conformal map.

In accordance with the definition of the space $\T$, it intrinsic {\bf Teichm\"{u}ller metric} on $\T$
is given by
$$
\tau_\T (\phi_\T (\mu), \phi_\T (\nu)) = \frac{1}{2} \inf \bigl\{
\tanh^{-1} \bigl( w^{\mu_*} \circ \bigl(w^{\nu_*} \bigr)^{-1} \bigr) : \ \mu_* \in \phi_\T(\mu), \nu_* \in \phi_\T(\nu) \bigr\},
$$
where $\phi_\T$ is the factorizing holomorphic projection $\Belt(D)_1 \to \T$. This metric is the
integral form of a canonically determined infinitesimal Finsler metric on the tangent bundle of $\T$.

Denote the equivalence classes (the points of $\T$) by $[w^\mu]$ and consider the balls
$$
B_\kappa(\mathbf 0, \kappa) = \{[w^\mu] \in \T: \ \tau_\T(\mathbf 0, [w^\mu]) < \kappa\},
\quad 0 < \kappa < \iy,
$$
centered at the origin $\mathbf 0 = \phi_\T(0)$.

As a complex manifold, the space $\T$ admits the hyperbolic invariant metrics. The largest among those is
the Kobayashi metric, and by the fundamental Royden-Gardiner theorem, this metric coincides with $\tau_\T$.
So the Teichm\"{u}ller metric arising from quasiconformality is also determined by appropriate holomorphic maps $\D \to \T$ (see \cite{EKK}, \cite{GL}, \cite{Le}, \cite{Ro}).

\bigskip
Using the {\bf Schwarzian derivatives}
$$
S_f(z) = \Bigl(\frac{f^{\prime\prime}(z)}{f^\prime(z)}\Bigr)^\prime - \frac{1}{2} \Bigl(\frac{f^{\prime\prime}(z)}{f^\prime(z)}\Bigr)^2
$$
of locally univalent holomorphic functions $f(z)$ on $\D$, which satisfy
 \be\label{1}
S_{f_1 \circ f}(z) = (S_{f_1} \circ f)(z) f^\prime(z)^2 + S_f(z)
\end{equation}
and belong to the complex Banach space $\B(\D)$ of hyperbolically bounded holomorphic functions (quadratic differentials) $\vp$ on $\D$ with norm
$$
\|\vp\|_\B = \sup_\D(1 - |z|^2)^2 |\vp(z)|,
$$
one obtains by applying the above factorization
$$
\mu \mapsto = S_{w^\mu}: \ \Belt(\D)_1 \to \T
$$
a holomorpic map of the ball $\Belt(\D)$ onto a bounded domain in the space $\B$ containing its origin,
which descends to a byholomorphic isomorphism between the space $\T$ and this domain.
It was established first by  Bers and gives a powerful tool in a lot of applications of Teichm\"{u}ller space theory.

\bigskip\noindent
{\bf 3. The necessary and sufficient conditions for biunivalence}.
The well-known result of Elisha Netanyahu \cite{Ne} states that
$$
\max_{\mathcal B}|a_2| = \fc{4}{3}.
$$
This estimate provides a {\bf necessary condition} for biunivalence.

A {\bf sufficient condition} is given by a result from \cite{Kr6}, which yields that univalent functions
$f$ with the second coefficient $|a_2|$ satisfying $|a_2| \le 1/2$ belong to $\mathcal B$.

More precisely, the following  covering lemma of Koebe's type is valid.
It concerns the (nonconstant) holomorphic maps $\chi$ from subdomains $G$ in complex Banach spaces $X = \{\mathbf x\}$
into the universal Teichm\"{u}ller space $\T$ modeled as a bounded subdomain of $\B$.
Suppose also that the image set $\chi(G)$ admits the circular symmetry, which means that for  every point $\vp \in \chi(G)$ the circle $e^{i \theta} \vp$ belongs entirely to this set.

Consider in the disk $\D$ the corresponding Schwarzian differential equations
 \be\label{2}
S_w(z) = \chi(\x)
\end{equation}
and pick their univalent in the entire disk solutions $w(z)$ satisfying $w(0) = 0, \ w^\prime(0) = 1$ (hence,  $w(z) = z  + \sum_2^\iy a_n z^n$). Put
$$
|a_2^0| = \sup \{|a_2|: \ S_w \in \chi(G)\},
$$
and let $w_0(z) = z + a_2^0 z^2 + \dots$ be one of the maximizing functions; its existence follows from compactness of the family of these $w(z)$ in topology of locally uniform convergence in $\D$ (note that
$a_2^0 \ne 0$). Then we have

\bigskip\noindent
{\bf Lemma 1}. \cite{Kr5} {\it For every indicated solution $w(z) = z + a_2 z^2 + \dots$ of the differential equation (2), the image domain $w(\D)$ covers entirely the disk $\{|w| < 1/(2 |a_2^0|)\}$.
The radius value $1/(2 |a_2^0|)$ is sharp for this collection of functions, and the circle $\{|w| = 1/(2|a_2^0|)\}$ contains points not belonging to $w(\D)$ if and only if $|a_2| = |a_2^0|$
(i.e., when $w$ is one of the maximizing functions). }

\bigskip
This implies

\bigskip\noindent
{\bf Lemma 2}. {\it All such univalent functions $f$ with $|a_2| \le 1/2$ belong to $\mathcal B$.
}

\bigskip
Lemma 1 yields also that the inverted functions
$$
W_w(\z) = 1/w(1/\z) = \z - a_2 + b_1 \z^{-1} + b_2 \z^{-2} + \dots, \quad |\z| > 1,
$$
map the complementary disk $\D^*$ onto the domains whose boundaries are entirely contained in the disk
$\{|W + a_2| \le \a_2^0|\}$.

Both estimates $|a_2| \le 4/3$ and $|a_2| \le 1/2$ are sharp; hence, these conditions cannot be improved
(only can be replaced by some equivalent conditions).

\bigskip
A standard criterion giving simultaneously the necessary and sufficient conditions for univalence  of a holomorphic function is obtained from the Grunsky inequalities arising from expansion
 \be\label{3}
\log \fc{f(z) - f(\z)}{z - \z} =  \sum\limits_{m, n = 1}^\iy c_{m n} z^m \z^n, \quad (z, \z) \in \D^2,
\end{equation}
where the principal branch of the logarithmic function is chosen. These coefficients satisfy
the inequality
 \be\label{4}
\Big\vert \sum\limits_{m,n=1}^\iy \ \sqrt{m n} \ c_{m n}(f) x_m x_n \Big\vert \le 1
\end{equation}
for any sequence $\mathbf x = (x_n)$ from  the unit sphere $S(l^2)$ of the Hilbert space $l^2$
with norm $\|\x\| = \bigl(\sum\limits_1^\iy |x_n|^2\bigr)^{1/2}$; conversely, this inequality
also is sufficient for univalence of a locally univalent function in $\D$ (see \cite{Gr}), \cite{Mi}, \cite{Po}).

It follows from (4) that $|c_{m n}| \le 1/\sqrt{m n}$ for all $m, n$.

Note that the inverted functions
$$
F_f(z) = 1/f(1/z) = z + b_0 + b_1 z^{-1} + \dots
$$
are $\hC$-holomorphic and univalent in the complementary disk $\D^* = \{z \in \hC: \ |z| > 1\}$ and have
the same Grunsky coefficients $c_{m n}$.

Generically, the inverse function $f^{-1}(w)$ is univalent only in some disk $\{|w| < r\}$ with $r < 1$,
and instead of (3),
$$
\log \fc{f^{-1}(w) - f^{-1}(\om)}{w - \om} = \sum\limits_{m, n = 1}^\iy c_{m n} r^{m+n} w^m \om^n.
$$
Accordingly, its univalence is equivalent to
$$
\sup \Big\{ \Big\vert \sum\limits_{m,n=1}^\iy \sqrt{m n} \ c_{m n} r^{m+n} x_m x_n \Big\vert : \
\mathbf x = (x_n) \in S(l^2) \Big\} \le 1.
$$
Hence, biunivalence is equivalent to having the same inequality (4) for both functions $f$ and $f^{-1}$, which also causes its rigidity.

\bigskip\noindent
{\bf 4. Main theorem}. Our first main result reveals a deep connection between biunivalence and geometry
of Teichm\"{u}ller balls.

\bigskip\noindent
{\bf Theorem 1}. {\it Any holomorphic function $f(z)$ on the disk $\D$, whose Schwarzian $S_f$ lies in
the ball $B_\T(\mathbf 0, \kappa)$ of radius $\kappa \le 1/4$ admits $\kappa$-quasiconformal extension
onto $\hC$ and is of the form
$$
f(z) = \gamma \circ f_0(z),
$$
where $f_0$ is biunivalent on $\D$ and $\gamma$ is a Moebius transformation of $\hC$.
The upper bound $1/4$ is sharp (cannot be increased).  }

\bigskip\noindent
{\bf Proof}. We represent the space $\T$ as a bounded domain in the space $\B$ of Schwazian derivatives
on $\D$ indicated in Section 2. Than any its point $\vp$ is the Schwarzian $\vp = S_f$ of a univalent function $f(z) = z + a_2 z^2 + \dots \in S$ admitting $k$-quasiconformal extension with $k < 1$.

We distinguish the functions normalized additionally by $f(\iy) = \iy$. For such functions, we have
the following

\bigskip\noindent
{\bf Lemma 3}. {\it For all functions $f(z) \in S$, admitting $k$-quasiconformal extension to $\hC$ with $f(\iy) = \iy$, we have the sharp bound
\be\label{5}
|a_2| \le 2 k,
\end{equation}
with equality only for the function }
$$
f_k(z) = \fc{z}{(1 - t k z)^2}, \quad |z| \le 1, \ \ |t| = 1.
$$

The proof of this estimate is given, for example, in \cite{KK}, Part 1.
Note that the extremal Beltrami coefficient among quasiconformal extensions of $f_k$ across the unit circle to $\hC$ is of the form
$$
\mu_0(z) = k |z|^3/z^3, \quad |z| > 1,
$$
so $|\mu_0| = \|\mu_0\|_\iy = \kappa$. Note also that by Ahlfors-Weill \cite{AW}, any holomorphic function $f(z)$ in the disk $\D$ satisfying
$$
(1 - |z|^2)^2 \ |S_f(z)| \le 2 k \quad \text{for all} z \in \D,
$$
with some  $k < 1$ is univalent on this disk and admits $k$-quasiconformal extension across $\mathbb S^1$ onto the whole sphere $\hC$.

(This gives only upper not sharp bound for the dilatation, because generically such quasiconformal extensions  of $f$ are not extremal, and the dilatation of the extremal Teichm\"{u}ller extension is strictly less than $k$.)

By construction of universal Teichm\"{u}ller space, all points $\vp$ of any ball $B_\kappa$ are the Schwarzian derivatives $S_f$ of univalent functions $f \in S$ admitting $\kappa$-quasiconformal extensions to $\hC$.
These functions being the solutions on $\D$ of the Schwarzian differential equations
$$
(w^{\prime\prime}/w^\prime)^\prime - \fc{1}{2} (w^{\prime\prime}/w^\prime)^2 = \vp,
$$
can be chosen so that their quasiconformal extensions satisfy  $f(\iy) = \iy$.

Taking $\kappa <1/2$, one obtains by Lemma 3 that the second coefficient $a_2(f)$ of any such function $f$ satisfies
$$
|a_2| \le 2 \kappa \le 1/2.
$$
By Lemma 2, this estimate implies that all such $f$ must cover the unit disk, and hence, their inverse functions are univalent in this disk.
This implies the assertion of the theorem.

\bigskip\noindent
{\bf 5. Other biholomorphic embeddings of universal Teichm\"{u}ller space}.
Now we dwell upon two other models of the space $\T$.

The first one is going back to J. Becker \cite{Be} and represents this space via bounded domain
$\mathbf b(\T)$ in the Banach space $\B_1(\D^*)$ of holomorphic functions $\psi$ on $\D^*$ with norm
$$
\|\psi\|= \sup_{\D^*} \ (|z|^2 -1 ) |z \psi(z)|.
$$
This domain is filled by functions $\psi_F = F^{\prime\prime}/F^\prime$ for univalent functions
$F(z) = z + b_0 + b_1 z^{-1} + \dots$ on $\D^*$ admitting quasiconformal extensions across $\mathbb S^1$.
In this case, one can use the inequality
$$
\Big\vert (|z|^2 - 1) \ z \fc{F^{\prime\prime}(z)}{F^\prime(z)} \Big\vert \le k, \quad z \in \D^*;
$$
with $k < 1$, which provides another sufficient condition for a holomorphic $f$ in $\D^*$ to be univalent and $k$-quasiconformally extended to $\hC$ (see \cite{Be}, \cite{Kr1}).

\bigskip
The third model of $\T$ was constructed by the author in \cite{Kr4}; it is determined by Grunsky coefficients of univalent functions.
In this model, the space $\T$ is represented as a bounded domain $\wt \T$ in some subspace of $l_\iy$.

Denote by $S_Q$ the collection of $f \in S$ admitting quasiconformal extensions onto $\D^*$
The Grunsky coefficients of $f \in S_Q$ span a $\C$-linear space $\mathcal L^0$ of sequences $\mathbf c = (c_{m n})$, which satisfy the symmetry relation $c_{m n} = c_{n m}$ and
$$
|c_{m n}| \le C(\mathbf c)/\sqrt{m n}, \quad C(\mathbf c) = \const < \iy \quad \text{for all} \ \ m, n \ge 1,
$$
with finite norm
 \be\label{6}
\|\mathbf c\| = \sup_{m,n} \ \sqrt{m n} \ |c_{m n}| +  \sup_{\mathbf x = (x_n) \in S(l^2)} \
\Big\vert \sum\limits_{m,n=1}^\infty \ \sqrt{mn} \ c_{mn} x_m x_n \Big\vert.
\end{equation}
Denote the closure of span $\mathcal L^0$ by $\mathcal L$ and note that the limits of convergent sequences
$\{\mathbf c^{(p)} = (c_{m n}^{(p)})\} \subset \mathcal L$ in the norm (6) also generate the double series
$\sum\limits_{m,n=1}^\infty \ c_{m n} z^m \z^n$
convergent absolutely in the unit bidisk $\{(z, \z) \in \C^2: \ |z| < 1, \ |\z| < 1\}$.

It was established in \cite{Kr4} that the sequences $\mathbf c$  corresponding to functions
$f \in S_Q$ fill a bounded domain $\wt \T$ in the indicated Banach space $\mathcal L$ containing the origin, and the correspondence
$$
S_f \longleftrightarrow \mathbf c(f) = (c_{m n}(f))
$$
determines a biholomorphism of this domain $\wt \T$ onto the space $\T$.

All three modeling domains $\T, \mathbf b(\T)$ and $\wt \T$ are biholomorphically equivalent, and their
basic geometric properties are similar.

In view of invariance of Teichm\"{u}ller metric under biholomorphic maps of $\T$, the above
theorem implies that also in both other domains $\mathbf b(\T)$ and $\wt T$ modeling $\T$ the hyperbolic ball of radius $1/4$ centered at the origin is filled by the points representing biunivalent functions on $\D$.

\bigskip\noindent
{\bf 6. Applications}.

\bigskip\noindent
{\bf 6.1. Covering theorems for univalent functions with quasiconformal extension}.
The class $S_k(\iy)$ of univalent functions in $\D$ having $k$-quasiconformal extension to $\hC$
with $f(\iy) = \iy$ is one of the basic classes in the distortion theory for univalent functions with quasiconformal extension.

We mention the following consequence of Theorem 1 giving a new explicit covering result.

{\bigskip\noindent}
{\bf Theorem 2}. {\it All functions $f \in S_k(\iy)$ with $k \le 1/4$ cover the unit disk}.

\bigskip
Another important class here is $\Sigma_k(0)$ of univalent functions $F(z) = z + b_0 + b_1 z^{- 1}+ \dots$
with quasiconformal extensions onto $\D$, satisfying $F(0) = 0$.

Using the translations $z \mapsto z + b$, one obtains the frames for the bounds of domains $F(\D^*)$ in the spirit of Lemma 2.

The covering theorems traditionally fill an important place in geometric function theory, starting from  the famous Koebe one-quoter theorem.

Theorem 2 involves a complete normalization of functions and is rather rigid.
More general covering theorems  for $k$-quasiconformally extendable univalent functions $f(z) = z + a_2 z^2 + \dots$ on the disk without an additional normalization have been obtained by Gutlyansky \cite{Gu1}, \cite{Gu2}, K\"{u}hnau \cite{Ku}, Schiffer-Schober \cite{SS1}, \cite{SS2}. These theorems also are based on estimating $a_2$.

\bigskip\noindent
{\bf 6.2. The distortion theory}. Theorem 1 also opens way to complete solving the general extremal
problems on classes $S_k(\iy)$, for $k \le 1/4$. Earlier only the special results have been obtained in
this direction.

The variational technique from \cite{Bel} and \cite{Kr1} gives the extremal functions in terms of
Beltrami coefficients of their inverses. On the other hand, an important consequence of biunvalence is
that the functions
$$
g(z) = z + \sum\limits_2^\iy b_n z^n,
$$
with the same coefficients $b_2, \ b_3, \dots$ as the inverse $f^{-1}$, also belong to $S_k(\iy)$.

Explicitly,
$$
b_2 = - a_2, \quad b_3 = 2 a_2^2 - a_3, \quad b_4 = - 5 a_2^3 + 5 a_2 a_3 - a_4, \dots
$$
This easily follows from the equality
$$
f^{-1}(w) = \fc{1}{2 \pi i} \int\limits_{|z| = 1 - \ve} \ \fc{z f^\prime(z)}{f(z) - w} dz,
$$
with arbitrary small $\ve > 0$.

We consider the general real or complex coefficient functionals
 \be\label{7}
J(f) = J(a_{n_1}, a_{n_2}, \dots, a_{n_s}),
\end{equation}
on $S_k(\iy)$, where $J$ is a continuously differentiable function of its arguments with $\grad J(f) \ne 0$,
$$
2 < n_1 < n_2, \dots, n_s < \iy.
$$
Noting that the indicated coefficients $b_n$ of $f^{-1}$ are polynomials of the initial coefficients $a_j$ of $f$, and vice versa, one obtains that this maximization problem is reduced to maximization of a functional $\wt J(f)$ on the same class $S_k(\iy)$ depending on the corresponding coefficients $b_n$, and  these functionals satisfy
$$
\max_{S_k(\iy)} |\wt J(f)| = \max_{S_k(\iy)} |J(g)|.
$$

Applying the variational method from \cite{Kr6}, one obtains in the same fashion as in \cite{Kr6} the following general theorem, which provides explicitly the Beltrami coefficients of all extremals of $J$:

\bigskip\noindent
{\bf Theorem 3}. {\it For any functional $J(f)$ of the form (7) and any $k \le 1/4$, we have the equalities
 \be\label{8}
\max_{S_k(\iy)} |J(f)| = \max_{S_k(\iy)} |\wt J(f)| = |\wt J(f^{\mu_k})|,
\end{equation}
where $\mu_k(z) = k |\vp_0(z)|/\vp_0(z), \quad z \in \D^*$ (extended by zero to $\D$),  with
 \be\label{9}
\vp_0(z) = \sum\limits_{l=1}^s \fc{\partial \wt J}{\partial b_{n_l}}(z).
\end{equation}
}

For example, in the case of $J(f) = a_3$, one must apply (8), (9) to $\wt J(f) = 2 b_2^2 - b_3$.

\bigskip\noindent
{\bf 7. Additional remarks}.

Together with Theorem 1, setting up biunivalence of functions $f \in S_k(\iy)$ with $k \le 1/4$,
Theorem 3 essentially strengthens many known distortion results for quasiconformally extentable univalent functions preserving the infinite point. Earlier such distortion theorems have been established only for sufficiently small dilatations $k \le \ve_0$ with implicit bound $\ve_0$.

Theorem 3 admits a straightforward extension to more general functionals 
$$
J(f) = J(a_{n_1}, a_{n_2}, \dots, a_{n_s}; f(z_1), f^\prime(z_1), \dots, f^{\a_1}(z_1); \dots, 
f(z_m), f^\prime(z_m), \dots, f^{\a_m}(z_1)), 
$$
where $z_1, \dots, z_m$ are the distinguished fixed points in $\D \setminus \{0\}$ with prescribed 
nonnegative orders $\a_1, \dots, \a_m$.

\bigskip
{\small\em{ \leftline{Department of Mathematics, Bar-Ilan University, 5290002 Ramat-Gan, Israel}
\leftline{and Department of Mathematics, University of Virginia,  Charlottesville, VA 22904-4137, USA}}

\end{document}